\documentclass[letterpaper, 10 pt, conference]{ieeeconf}  

\IEEEoverridecommandlockouts                              
\usepackage{amssymb,amsbsy,amsmath,amsfonts,amssymb,amscd,mathrsfs}

\newtheorem{theorem}{Theorem}[section]

\newtheorem{lemma}[theorem]{Lemma}

\newtheorem{proposition}[theorem]{Proposition}
\newtheorem{definition}[theorem]{Definition}

\def\vrexp{\stackrel{\longrightarrow}{\exp}}
\newcommand{\rexp}[3]{\,  \vrexp \int_0^{#1} #2 \, d #3}
\newcommand{\deinde}[1]{\frac{\partial}{\partial #1}}

\newcommand{\ad}{\operatorname{ad}}
\newcommand{\vep}{\varepsilon}
\newcommand{\eps}{\varepsilon}
\renewcommand{\H}{\mathcal{H}}

\newcommand{\pro}{\Upsilon} 

\newcommand{\spann}{\mathrm{span}}
\newcommand{\NN}{\mathbb{N}}
\newcommand{\CC}{{\mathbb{C}}}
\newcommand{\RR}{{\mathbb{R}}}

\renewcommand{\l}{\ell}

\title{\LARGE \bf
Controllability of the bilinear Schr\"odinger equation with several controls and application to a 3D molecule$^{*}$}

\author{Ugo Boscain$^{1}$, Marco Caponigro$^2$, and Mario Sigalotti$^3$
\thanks{$^{*}$ This research has been supported by the European Research Council, ERC
StG 2009 ``GeCoMethods'', contract number 239748, by the ANR project
GCM, program ``Blanche'',
project number NT09-504490}%
\thanks{$^1$ Ugo Boscain is with 
Centre National de Recherche Scientifique (CNRS), CMAP, \'Ecole Polytechnique, Route de Saclay, 91128 Palaiseau Cedex, France, and Team GECO, INRIA-Centre de Recherche Saclay
{\tt \small ugo.boscain@polytecnique.edu}}%
\thanks{$^2$ Marco Caponigro is with
Department of Mathematical Sciences and Center for Computational and Integrative Biology,
Rutgers - The State University of New Jersey, Camden NJ 08102, USA
{\tt \small marco.caponigro@rutgers.edu}}%
\thanks{$^3$ Mario Sigalotti is with
INRIA-Centre de Recherche Sacaly, Team GECO and 
CMAP, \'Ecole Polytechnique, Route de Saclay, 91128 Palaiseau Cedex, France
{\tt \small mario.sigalotti@inria.fr}}%
}

\begin{document}

\maketitle
\thispagestyle{empty}
\pagestyle{empty}

\begin{abstract}
We show the approximate rotational controllability of  a polar linear molecule by means of three nonresonant linear polarized laser fields.
The result is based on a
general approximate controllability result for the
bilinear Schr\"odinger equation, 
with wavefunction varying in the unit sphere of an infinite-dimensional Hilbert space and with several control potentials, under the assumption that the internal Hamiltonian has discrete spectrum. 
\end{abstract}


\section{Introduction}

Rotational molecular dynamics is one of the most important examples of quantum systems with an infinite-dimensional Hilbert space and a discrete spectrum. 
Molecular orientation and alignment are well-established topics in the quantum control of molecular dynamics both from the experimental and theoretical points of view (see \cite{seideman,stapelfeldt} and references therein). For linear molecules driven by linearly polarized laser fields in gas phase, alignment means an increased probability direction along the polarization axis whereas orientation requires in addition the same (or opposite) direction as the polarization vector. Such controls have a variety of applications extending from chemical reaction dynamics to surface processing, catalysis and nanoscale design. A large amount of numerical simulations have been done in this domain but the mathematical part is not yet fully understood. From this perspective, the controllability problem is 
a necessary step towards comprehension.

 We focus in this paper on the control by laser fields of the rotation of a rigid linear molecule in $\RR^3$. This control problem corresponds to the control of the Schr\"odinger equation on the unit sphere $S^2$. 
 We show that the system driven by three fields along the three axes is  
 approximately  controllable
 for arbitrarily small controls. 
 This means, in particular,  that there exist control strategies which bring the initial state arbitrarily close to states maximizing the molecular orientation \cite{sugny}.
\subsection{The model}
We consider a polar linear molecule in its ground vibronic state subject to three nonresonant (with respect to the vibronic frequencies) linearly polarized laser fields. 
The control is given by the electric fields $E=(u_1,u_2,u_3)$ depending on time and constant in space. We neglect in this model the polarizability tensor term which corresponds to the field-induced dipole moment. This approximation is correct if the intensity of the laser field is sufficiently weak. Despite its simplicity, this equation reproduces very well the experimental data on the rotational dynamics of rigid molecules (see \cite{stapelfeldt}).

Up to normalization of physical constants (in particular, in units such that $\hbar=1$), the dynamics is ruled by the equation
\begin{align}
i\frac{\partial\psi(\theta,\varphi,t)}{\partial t}=& -\Delta \psi(\theta,\varphi,t) + (u_1(t) \sin \theta\cos \varphi\nonumber\\
&+ u_2(t) \sin \theta\sin \varphi+u_3(t) \cos \theta)\psi(\theta,\varphi,t)
\label{1}
\end{align}
where $\theta,\varphi$ are the spherical coordinates, which are related to the Euclidean coordinates by the identities 
$$
x = \sin \theta \cos \varphi, \quad y = \sin \theta\sin \varphi,\quad z = \cos \theta,
$$
while 
$\Delta$ is the Laplace--Beltrami operator on the sphere (called in this context the \emph{angular momentum operator}), 
i.e., 
$$
\Delta = \frac{1}{\sin \theta} \deinde{\theta} \left( \sin \theta \deinde{\theta}\right) + \frac{1}{\sin^2 \theta} \frac{\partial^2}{\partial \varphi^2}.
$$

The wavefunction $\psi(\cdot,\cdot,t)$ evolves in the unit sphere ${\cal S}$ of ${\cal H}=L^2(\mathbb{S}^2,\CC)$.

\subsection{The main results}\label{SEC_main_result}
In the following we denote by $\psi(T;\psi_0,u)$ the solution at time $T$ of equation \eqref{1}, corresponding to control $u$ and   with initial condition $\psi(0;\psi_0,u)=\psi_0$, belonging to ${\cal S}$.

Our main result says that \eqref{1} is approximately controllable with arbitrarily small controls.
\begin{theorem}
\label{controllability}
For every  $\psi^0$, $\psi^1$  belonging to ${\cal S}$  and every $\eps,\delta_1,\delta_2,\delta_3>0$,
 there exist $T>0$ and $u\in L^\infty([0,T],[0,\delta_1]\times [0,\delta_2]\times [0,\delta_3])$ such that
 $\|\psi^1-\psi(T;\psi^0,u)\|<\eps$.
\end{theorem}

The proof of the result is based on arguments inspired by those developed in \cite{noi,ancoranoi}. 
There are two main difficulties preventing us
to apply those results to the case under consideration: firstly, we deal here with several control parameters, while those general results were specifically conceived for the single-input case. Notice that, because of symmetry obstructions, equation \eqref{1} is not controllable with only one of the three controls $u_1$, $u_2$, $u_3$.  
Secondly, the general theory developed in \cite{noi,ancoranoi} is based on 
nonresonance conditions on the spectrum of the drift Schr\"odinger operator (the \emph{internal Hamiltonian}). 
The Laplace--Belatrami operator on $S^2$, however, has a severely degenerate spectrum. It is known, indeed, that the $\ell$-th eigenvalue $-i\ell(\ell+1)$ has multiplicity $2\ell+1$.
In \cite{noi} we proposed a perturbation technique in order to overcome resonance relations in the spectrum of the drift. This technique was applied in \cite{noicdc} to the case of the orientation of a molecule confined in a plane driven by one control.  The planar case is already technically 
challenging and a generalization to the case of three controls in the space will hardly provide an apophantic proof of the approximate controllability result. 
We therefore provide a general multi-input result which can be applied to the control problem defined in \eqref{1}, up to the computation of certain Lie algebras associated with its Galerkin approximations.

The structure of the paper
is the following: in the next section
we present the general multi-input abstract framework and  we recall some previously known controllability and non-controllability results.
In Section~\ref{III} we 
prove our  main sufficient condition for approximate controllability. 
Finally, in Section~\ref{rotational} we prove that the abstract result applies to system \eqref{1}.

\section{Abstract framework}

\begin{definition}\label{def:system}
Let $\cal H$ be an infinite-dimensional Hilbert space with scalar product $\langle \cdot ,\cdot \rangle$ and  $A,B_{1}, \ldots,B_{p}$ be (possibly unbounded) linear operators on $\H$,
with domains $D(A), D(B_{1}), \ldots, D(B_{p})$. Let $U$ be a subset of $\mathbb{R}^{p}$.
Let us introduce the 
controlled equation
\begin{equation} \label{eq:main}
\frac{d\psi}{dt}(t)=(A+u_1(t)B_1+ \cdots + u_p(t) B_p) \psi(t),  u(t) \in U \subset \RR^{p}.
\end{equation}
We say that $(A,B_{1}, \ldots, B_{p},U,\Phi)$ satisfies $(\mathfrak{A})$ 
if the following assumptions are verified:
\begin{description}
\item[($\mathfrak{A}1$)] $ \Phi = (\phi_k)_{k \in \mathbf{N}}$ is an Hilbert basis of $\cal H$
made of eigenvectors of $A$ associated with the family of eigenvalues $(i \lambda_{k})_{k \in \mathbb{N}}$; 
\item[($\mathfrak{A}2$)] $\phi_k \in D(B_{j})$ for every $k \in
\mathbb{N}, j=1,\ldots,p $;
\item[($\mathfrak{A}3$)] 
$A+u_1B_1+ \cdots + u_p B_p:\mathrm{span}\{\phi_k \mid k\in\mathbb{N}\}\to \mathcal{H}$ is essentially skew-adjoint for every $u\in U$;
\item[($\mathfrak{A}4$)] if $j\neq k$ and $\lambda_{j} = \lambda_{k}$ then $\langle \phi_{j},B_{l}\phi_{k} \rangle = 0$ for every $l=1,\ldots,p$.
\end{description}
\end{definition}

If $(A,B_{1}, \ldots, B_{p},U,\Phi)$ satisfies $(\mathfrak{A})$ then, for every $(u_{1},\ldots,u_{p}) \in \mathbb{R}^{p}$, $A+u_{1}B_{1} + \cdots + u_{p} B_{p}$ generates a unitary group $e^{t(A+u_{1}B_{1} + \cdots + u_{p} B_{p})}$. It is therefore possible to define the propagator $\pro^{u}_{T}$ at time $T$ of system~\eqref{1} associated with a $p$-uple of piecewise constant controls $u(t)=(u_{1}(t), \ldots, u_{p}(t))$ by concatenation. If, moreover, the potentials $B_{1},\ldots, B_{p}$ are bounded operators then the definition can be extended by continuity to every $L^{\infty}$ control law.

\begin{definition}\label{def:controllability}
Let $(A,B,U,\Phi)$ satisfy $(\mathfrak{A})$.
We say that \eqref{eq:main} is \emph{approximately controllable}
if for every $\psi_0,\psi_1$ in the unit sphere of $\cal H$
 and every $\vep>0$ there exist a piecewise constant control function $u:[0,T] \to U$ such that 
$
\|\psi_1-  \pro^{u}_{T}(\psi_{0})\| <\vep.
$
\end{definition}

\begin{definition} \label{DEF_simultaneous_controllability}
 Let $(A,B,U,\Phi)$ satisfy $(\mathfrak{A})$ . 
We say that \eqref{eq:main} is \emph{approximately simultaneously controllable} if for every $r$ in $ \NN$, $\psi_1,\ldots,\psi_{r}$ in $\H$, $\hat\Upsilon$ in $\mathbf{U}(\H)$, and $\vep>0$  there exists a piecewise constant control $u:[0,T]\rightarrow U$
 such that
$$
\left \| \hat\Upsilon \psi_k - \pro^{u}_T \psi_k \right \|<\vep,\qquad k=1,\dots,r.
$$
\end{definition}

\subsection{Short review of controllability results}

The controllability of system~\eqref{eq:main} is a well-established topic when the state space $\mathcal{H}$ is 
finite-dimensional (see for instance \cite{dalessandro-book} and reference therein),
thanks to general controllability methods for left-invariant control systems on compact Lie groups (\cite{brock,jur}).

When $\mathcal{H}$ is infinite-dimensional, 
it is known 
that the bilinear Schr\"odinger equation is not controllable (see \cite{BMS,turinici}). Hence, one has to look for weaker controllability properties as, for instance, approximate controllability or controllability between  eigenstates of the Sch\"odinger operator (which are the most relevant physical states).
In certain cases where 
the dimension of the domain where the controlled PDE is defined is equal to one  a description of the reachable set has been provided~\cite{Beauchard1, beauchard-coron, camillo}.
For dimension larger than one 
or for more general situations, the exact description of the reachable set appears to be more difficult and at the moment only approximate controllability results are available. 
Most  of them are for the single-input case 
(see, in particular, \cite{beauchard-nersesyan, ancoranoi, noi, mirrahimi-aihp, Nersy, fratelli-nersesyan, nersesyan}), except for some 
  approximate controllability result for specific systems (\cite{ervedoza_puel}) and some general approximate controllability result 
  between eigenfunctions based on adiabatic methods \cite{adiabatiko}.

\subsection{Notation}

Set $b^{(l)}_{jk} = \langle \phi_{j}, B_{l} \phi_{k} \rangle$, $l=1,\ldots,p$.
For every $n$ in $\NN$, define the orthogonal projection
$$
 \pi_n: \H \ni \psi\mapsto \sum_{j\leq n} \langle \phi_j,\psi\rangle
\phi_j \in \H.
$$

Given a linear operator $Q$ on $\H$ we identify the linear operator
$\pi_{n} Q \pi_{n}$ preserving 
$\spann\{\phi_{1},\ldots, \phi_{n}\}$ with 
its  
$n \times n$ complex matrix representation with 
respect to the basis $(\phi_{1},\ldots, \phi_{n})$.

\section{Main abstract controllability result in the multi-input case}
\label{III}

Let us introduce the set $\Sigma_{N}$ of spectral gaps associated with the $N$-dimensional Galerkin approximation as
$$
\Sigma_{N} = \{|\lambda_j - \lambda_k|\,\mid \, j,k = 1, \ldots, N, \lambda_j \neq \lambda_k\}.
$$

For every $\sigma \in \Sigma_{N}$, let
$$
B^{(N)}_{\sigma} (v_1, \ldots, v_p)_{j,k}=  (v_1B^{(N)}_1+ \ldots + v_p B^{(N)}_p)_{j,k} \delta_{\sigma, |\lambda_j - \lambda_k|}.
$$
The $N\times N$ matrix $B^{(N)}_{\sigma} (v_1, \ldots, v_p)$ corresponds to the choice of the controls $v_1, \ldots, v_p$ and to the ``activation" of the spectral gap $\sigma$. 
Define 
$$
\mathcal{M}_N =  \{B^{(N)}_{\sigma} (v_1, \ldots, v_p)\mid \sigma \in \Sigma_{N}, v_1,\dots, v_p \in[0,1]\}
$$
and 
\begin{align*}
& \mathcal{M}^{n}_0 = 
\left\{A^{(n)} - \frac{\mathrm{tr}(A^{(n)})}{n}I_{n} \right\}
\cup \\
&\, \left\{M \in\mathfrak{su}(n) \mid \forall N \geq n\, \exists\,  Q \in  \mathcal{M}_N \mbox{ s.t. } Q = 
\left(
\begin{array}{c|c}
M&0\\ \hline 0 & *
\end{array}
\right)
\right\}.
\end{align*}
The set $\mathcal{M}^{n}_0$ represents 
``compatible dynamics" for the $n$-dimensional Galerkin approximation (\emph{compatible}, that is, with higher dimensional Galerkin approximations). 

\begin{theorem}[Abstract multi-input controllability result]\label{metalemma}
Let $U=[0,\delta]^p$ for some $\delta>0$. 
If for every $n_0\in \mathbb{N}$ there exist $n> n_0$ 
such that
\begin{equation}\label{hypothesis}
\mathrm{Lie} \mathcal{M}^{n}_0 = \mathfrak{su}(n),
\end{equation}
then the  system
$$
\dot{x} = (A+ u_1B_1+ \cdots + u_p B_p)x,\quad u\in U,
$$
is approximately simultaneously controllable.
\end{theorem}

\subsection{Preliminaries}

The following technical result, which we shall use in the proof of Theorem~\ref{metalemma}, has been proved in \cite{ancoranoi}. 

\begin{lemma}\label{lem:convexification}
Let $\kappa$ be a positive integer and $\gamma_{1}, \ldots, \gamma_{\kappa} \in \RR\setminus \{0\}$ be  such that
$|\gamma_{1}|\neq|\gamma_{j}|$ for $j =2,\ldots, \kappa.$
Let 
$$ 
\varphi(t) = (e^{it\gamma_{1}}, \ldots, e^{it\gamma_{\kappa}}).
$$
Then, for every $\tau_0\in \RR$, we have
$$ 
\overline{\mathrm{conv}{\varphi([\tau_0,\infty))}} \supseteq \nu \mathbb{S}^{1} \times \{(0, \ldots, 0)\}\,,
$$
where
$
\nu=\prod_{k=2}^{\infty} \cos \left (\frac{\pi}{2 k} \right ) >0.
$
Moreover,
for every $R>0$ and $\xi \in \mathbb{S}^{1}$ there exists a sequence $(t_{k})_{k\in \NN}$ such that
$
t_{k+1} - t_{k} > R
$
and 
$$
\lim_{h\to \infty } \frac{1}{h} \sum_{k=1}^{h} \varphi (t_{k}) = (\nu \xi,0,\ldots,0)\,.
$$
\end{lemma}

\subsection{Time reparametrization}

For every piecewise constant function $z(t) = \sum_{k=1}^{K} z_{k} \chi_{[s_{k-1},s_{k})} (t)$
such that
$z_{k}> 0$, for every $k=1,\ldots, K$,
 and $v_{j}(t) = \sum_{k=1}^{K} v^{(j)}_{k} \chi_{[s_{k-1},s_{k})} (t)$ with $j=1,\ldots,p$,
we consider the system
\begin{equation} \label{eq:main-repar}
\frac{d\psi}{dt}(t)=(z(t)A+v_1(t)B_1+ \cdots + v_p(t) B_p) \psi(t).
\end{equation}

System \eqref{eq:main-repar} can be seen as a time-reparametrisation of system~\eqref{eq:main}. 
Let $\psi(t)$ be the solution of \eqref{eq:main} with initial condition $\psi_{0}\in \H$ associated with the piecewise constant control $u(\cdot)$ with components  
$u_{j}(\cdot) = \sum_{k=1}^{K}u^{(j)}_{k}\chi_{[t_{k-1},t_{k})}(\cdot)$, $j=1,\ldots,p$. 
If $s_{k} = \frac{t_{k}- t_{k-1}}{z_{k}} + s_{k-1}$, $s_{0}=0$, $v^{(j)}_{k} =u^{(j)}_{k} z_{k} $ for every $k=1,\ldots,K$, $j=1,\ldots,p$, then the solution
$\tilde\psi(t)$ of \eqref{eq:main-repar} with the initial condition $\psi_{0}\in \H$ associated with the
controls $z(t), v_{1}(t), \ldots, v_{p}(t)$ satisfies
$$
\tilde\psi \left( \int_{0}^{t} \sum_{k=1}^{K} \frac{1}{z_{k}} \chi_{[t_{k-1},t_{k})}(s)ds\right)= \psi(t)\,.
$$
Controllability issues for system~\eqref{eq:main} and~\eqref{eq:main-repar} are equivalent. Indeed, consider piecewise constant controls
$z:[0,T_{v}] \to [1/\delta,\infty)$, $z(t)=\sum_{k=1}^{K} z_{k} \chi_{[s_{k-1},s_{k})} (t)$ and
$v_{j}:[0,T_{v}] \to [0,1]$,
$v_{j}(t) = \sum_{k=1}^{K} v^{(j)}_{k} \chi_{[s_{k-1},s_{k})} (t)$ with $j=1,\ldots,p$,
achieving controllability (steering system~\eqref{eq:main-repar} from $\psi_{j}$ to $\hat\Upsilon \psi_{j}$, $j=1,\ldots,r$ in a time $T_{v}$) . Then the controls $u_{j}(t) = \sum_{k=1}^{K}u^{(j)}_{k}\chi_{[t_{k-1},t_{k})}$, $j=1,\ldots,p$ defined by $u_{k}^{(j)}=v_{k}^{(j)}/z_{k}$ and $t_{0}=0, t_{k}=(s_{k}-s_{k-1})z_{k}+t_{k-1}$, steer system~\eqref{eq:main} from $\psi_{j}$ to $\hat\Upsilon \psi_{j}$, $j=1,\ldots,r$ in a time $T_{u}$.

\subsection{Interaction framework}

Let $\omega(t) = \int_{0}^{t}z(s)ds$, and 
$w_{j}(t) = \int_{0}^{t} v_{j}(s)ds$ for $j=1,\ldots,p$. 
Let $\psi(t)$ be the solution of \eqref{eq:main-repar} with  initial condition $\psi_{0}\in \H$ associated with the
controls $z(t), v_{1}(t), \ldots, v_{p}(t)$ and set
$$
y(t) = e^{-\omega(t)A } \psi(t).
$$
For $\omega, v_{1},\ldots,v_{p} \in \RR$ set $ \Theta(\omega, v_{1},\ldots,v_{p})=e^{-\omega A  }
(v_{1} {B_{1}} + \cdots + v_{p}(t) {B_{p}} )
e^{\omega A }$,
then $y(t)$ satisfies
\begin{equation}\label{eq:main-interaction}
\dot y(t)   = 
\Theta(\omega(t), v_{1}(t),\ldots,v_{p}(t))y(t).
\end{equation}

Note that
\begin{align*}
\Theta(\omega, v_{1},\ldots,v_{p})_{jk} &= \langle \phi_{k}, \Theta(\omega, v_{1},\ldots,v_{p}) \phi_{j} \rangle\\
&=
e^{i(\lambda_{k} - \lambda_{j})\omega} \left(v_{1} b_{jk}^{(1)}+\cdots + v_{p} b_{jk}^{(p)} \right).
\end{align*}

Notice that
$
|y(t)| = |\psi(t)|,
$
for every $t \in [0,T_{v}]$ and for every $(p+1)$-uple of piecewise constant controls 
$z : [0,T_{v}] \to [1/\delta,+\infty)$,
$v_{1}, \ldots, v_{p}:[0,T_{v}] \to [0,1]$.

\subsection{Galerkin approximation}

\begin{definition}
Let $N \in \NN$.  The \emph{Galerkin approximation}  of \eqref{eq:main-interaction}
of order $N$ is the system in $\H$
\begin{equation}\label{eq:galerkin-N}
\dot x = \Theta^{(N)}(\omega, v_{1},\ldots,v_{p}) x 
\end{equation}
where $\Theta^{(N)}(\omega, v_{1},\ldots,v_{p})=\pi_N \Theta(\omega, v_{1},\ldots,v_{p}) \pi_N = \left( 
\Theta(\omega, v_{1},\ldots,v_{p})_{jk} \right)_{j,k =1}^{N}$.
\end{definition}

\subsection{First step: choice of the order of the Galerkin approximation
}\label{sec:n}

In order to prove approximate simultaneous controllability,  
we should take $r$ in $ \NN$, $\psi_1,\ldots,\psi_{r}$ in $\H$, $\hat\Upsilon$ in $\mathbf{U}(\H)$, and $\vep>0$ and prove 
 the existence of a piecewise constant control $u:[0,T]\rightarrow U$
 such that
$$
\left \| \hat\Upsilon \psi_k - \pro^{u}_T \psi_k \right \|<\vep,\qquad k=1,\dots,r.
$$

Notice that for  $n_{0}$ large enough there 
exists  $U \in SU(n_{0})$ such that 
$$
| \langle \phi_{j},\hat\Upsilon \psi_k \rangle  -  \langle \pi_{n_{0}}\phi_{j},U \pi_{n_{0}}\psi_k \rangle| <\vep
$$ 
for every  $1 \leq k\leq r$ and $j\in\NN$.
This simple fact 
suggest to prove 
approximate simultaneous controllability
by studying the controllability of 
\eqref{eq:galerkin-N} 
in the Lie group $SU(n_{0})$.

\subsection{Second step: control in $SU(n)$}

Let $n\geq n_0$ satisfy hypothesis \eqref{hypothesis}. It follows from standard controllability results on compact Lie groups (see \cite{jur}) that  
for every $U \in SU(n)$ 
there exists a path $M:[0,T_{v}] \to \mathcal{M}^{n}_0$ such that 
$$
 \rexp{T_{v}}{M(s)}{s} = U,
$$
where the chronological notation $\rexp{t}{V_s}{s}$ is used for the flow from time $0$ to $t$ of the time-varying equation $\dot q=V_s(q)$ (see \cite{book2}). More precisely, there exists a finite partition in intervals $(I_{k})_{k}$ of $[0,T_{v}]$ such that
for every  $t \in I_{k}$
either there exist  $v_{1}, \ldots, v_{p}\in [0,1]$ and $\sigma \in \Sigma_{N}$ such that
$$
 M(t) = \pi_{n} B^{(N)}_{\sigma}(v_{1}, \ldots, v_{p}) \pi_{n},
$$
or 
$$
M(t) =  A^{(n)} - \frac{\mathrm{tr}(A^{(n)})}{n}I_{n}.
$$
In particular, 
\begin{equation}\label{eq:BN}
M(t)_{j,k} = 0, \quad \mbox{ for every } t \in [0,T_{v}],  j \leq n,  k>n.
\end{equation}

\subsection{Third step: control of $\mathcal{M}_{N}$}

\begin{lemma}
For every $N \in \NN$, $\delta >0$, and for every piecewise constant $v_{1}, \ldots, v_{p}: [0,T_{v}] \to [0,1]$ and $\sigma:[0,T_{v}] \to \Sigma_{N}$ there exists a sequence
$(z_{h}(\cdot))_{h\in \NN}$ of
 piecewise constant functions from $[0,T_{v}]$ to $[1/\delta, \infty)$, such that
\begin{align*}
 &\left\|\int_{0}^{t} \Theta^{(N)}(z_{h}(s),v_{1},\ldots,v_{p}) ds \right.\\
& \quad \left. - \int_{0}^{t}B^{(N)}_{\sigma(s)}(v_{1}(s),\ldots,v_{p}(s))ds\right\| \to 0
\end{align*}
uniformly with respect to $t \in [0,T_{v}]$ as $h$ tends to infinity. 
\end{lemma}

In other words, every piecewise constant path in $\mathcal{M}_{N}$ can be approximately tracked by
system~\eqref{eq:galerkin-N}.

\noindent{\bf Proof.}
Fix $N \in \NN$.
We are going to construct the control  
$z_{h}$ by applying recursively Lemma~\ref{lem:convexification}. Consider an interval $[t_{k},t_{k+1})$ in which 
$v_{j}(t)$, $j=1,\ldots,p$, and $\sigma(t)$ are constantly equal to $v_{j} \in [0,1]$, $j=1,\ldots,p,$ and $\sigma \in \Sigma_{N}$ respectively.
Apply Lemma~\ref{lem:convexification} with $\gamma_{1} = \sigma$, 
$\{\gamma_{2}, \ldots,\gamma_{\kappa}\} = 
\Sigma_{N} \setminus \{\sigma\}$,
$R=T$
and $\tau_0 = \tau_{0}(k)$ to be fixed later depending on $k$. Then, for every $\eta >0$, there exist 
$h = h(k) > 1/\eta$ and a sequence $(w_{\alpha}^{k})_{\alpha=1}^{h}$ such that $w_{1}^{k} \geq t_{0}$, $w_{\alpha}^{k}-w_{\alpha-1}^{k} > R$, and such that 
\begin{align*}
&\left|\frac{1}{h} \sum_{\alpha=1}^{h} e^{i(\lambda_{l} - \lambda_{m}) w^{k}_{\alpha} } \right.  \\
& \quad \left. - \nu \frac{(v_{1}\overline{B_{1}}^{(N)} + \ldots + v_{p}\overline{B_{p}}^{(N)})_{l,m}}{|(v_{1}B_{1}^{(N)} + \ldots + v_{p}B_{p}^{(N)})_{l,m}|} \delta_{\sigma, |\lambda_{l} - \lambda_{m}|}
\right| < \eta,
\end{align*}

Set $\tau^{k}_{\alpha} = t_{k}+ (t_{k+1} - t_{k}){\alpha/h}$, $\alpha = 0,\ldots, h$, and define the piecewise constant function
\begin{equation}\label{eq:vht}
\omega_{h}(t) = \sum_{k\geq 0} \sum_{\alpha=1}^{h(k)} 
w_{\alpha}^{k}  \chi_{[\tau^{k}_{\alpha-1}, \tau^{k}_{\alpha})}(t)\,.
\end{equation}
 Note that by choosing $\tau_{0}(k) = w^{k-1}_{h(k-1)} + R$ for $k\geq 1$ and $\tau_{0}(0) = R$ 
we have that $\omega_{h}(t)$ is non-decreasing.

Following the smoothing procedure of~\cite[Proposition~5.5]{ancoranoi} one can construct  the desired sequence of control $z_{h}$. The idea is to approximate $\omega_{h}(t)$ by suitable piecewise linear functions with slope greater than $1/\delta$. Then $z_{h}$ can be constructed from the derivatives of these functions. 
\hfill$\Box$

As a consequence of last proposition by~\cite[Lemma~8.2]{book2} we have that
\begin{align*}
& \left\| 
 \rexp{t}{\Theta^{(N)}(z_{h}(s),v_{1}(s),\ldots,v_{p}(s)}{s}\right. \\ 
& \quad -\left. \rexp{t}{B^{(N)}_{\sigma(s)}(v_{1}(s),\ldots,v_{p}(s))}{s}
 \right\|
 \to 0
\end{align*}
uniformly with respect to $t \in [0,T_{v}]$  as $h$ tends to infinity.

\subsection{Fourth step: control of the infinite-dimensional system}
Next proposition states that, roughly speaking, we can pass to the limit as $N$ tends to infinity without losing the controllability property proved for the finite-dimensional case. Its proof can be found in~\cite[Proposition~5.6]{ancoranoi}. It is based on the particular form~\eqref{eq:BN} of the operators involved,  
since the fact that the operator has several zero elements guarantees that the difference between the dynamics of the infinite-dimensional system and the dynamics of the Galerkin approximations is small.

\begin{proposition}
For every $\vep>0$, for every $\delta>0$, and for every trajectory $U \in SU(n)$  there exist  
piecewise constant controls $u_{j}: [0,T_{u}]\to [0,\delta], j=1,\ldots,p$ such that 
 the associated propagator $\Upsilon^{u}$ of \eqref{eq:main} satisfies 
$$
\big|  | \langle \pi_{n} \phi_{j} , U\pi_{n}\phi \rangle | - | \langle \phi_j , \Upsilon^{u}_{T_{u}} \phi\rangle | \big| <\vep
$$
 for every $\phi\in\spann \{\phi_1,\dots,\phi_n\}$ with $\|\phi\|=1$ and for every  $j$ in $\NN$.
\end{proposition}

We recall now a controllability result for the phases (see~\cite[Proposition 6.1 and Remark 6.3]{ancoranoi}).
This property, stated in the proposition below, together with the controllability up to phases
proved in the previous section,  is sufficient to conclude the proof of Theorem~\ref{metalemma}.

\begin{proposition}
Assume that, for every $\hat\Upsilon \in \mathbf{U}(\H)$, $m $ in $\mathbf{N}$, $\delta >0$, and $\vep>0$, there exist $T_{u}>0$ and  piecewise constant controls $u_{j}:[0,T_u]\rightarrow [0,\delta]$, $j=1,\ldots,p$ such that
the associated propagator $\pro^{u}$ of equation~\eqref{eq:main}
satisfies
$$
\big |  |\langle   \phi_{j} ,\hat\Upsilon \phi\rangle | - |\langle   \phi_{j},\pro^{u}_{T_{u}}\phi \rangle |\big |< \vep, 
$$
for every $j \in \mathbf{N}$
and 
$\phi\in\mathrm{span}\{\phi_1,\dots,\phi_m\}$ with $\|\phi\|=1$.
Then~\eqref{eq:main} is simultaneously approximately controllable. 
\end{proposition}

\section{$3$D molecule}
\label{rotational}

Let us go back to the system presented in the introduction for the orientation of a linear molecule,   
\begin{equation}
i \hbar \dot \psi = -\Delta \psi + (u_1 \cos \theta + u_2 \cos \varphi\sin \theta+ u_3 \sin \varphi\sin \theta)\psi,
\end{equation}
where $\psi(t) \in \H=L^2(\mathbb{S}^2,\CC)$.

A basis of eigenvectors of the Laplace--Beltrami operator $\Delta$ is given by  
the spherical harmonics $Y^m_\l(\theta, \varphi)$, which sastisfy
$$
\Delta Y^m_\l (\theta,\varphi) = -\l(\l+1)Y^m_\l(\theta, \varphi).
$$

We are first going to prove that for every $\l\in\NN$ the system projected on the $(4\l+4)$-dimensional linear  space
$$
\mathcal{L}:=\mathrm{span}\{Y^{-\l}_\l,\ldots,Y^\l_\l,Y^{-\l-1}_{\l+1},\ldots,Y^{\l+1}_{\l+1}\}
$$
is controllable. 
More precisely, chosen a reordering $(\phi_k)_{k\in \NN}$ of the spherical harmonics in such a way that 
$$\{\phi_k\mid k=1,\dots,4\l+4\}=\{Y^{-\l}_\l,\ldots,Y^\l_\l,Y^{-\l-1}_{\l+1},\ldots,Y^{\l+1}_{\l+1}\},$$
we are going to prove that 
$$
\mathrm{Lie} \mathcal{M}^{4\l+4}_0 = \mathfrak{su}(4\l+4).
$$

\subsection{Matrix representations}

Denote
by $J_\l$ 
the set of integer pairs 
$\{(j,k)\mid j=\l,\l+1,\ k=-j,\dots,j\}$. 
Consider an ordering $\omega:\{1,\dots,4\l+4\}\to J_\l$. 
Let $e_{j,k}$ be the $(4\l+4)$-square matrix whose entries are all zero, but the one at line $j$ and column $k$ 
which is equal to $1$.
Define
$$
E_{j,k} = e_{j,k} - e_{k,j}, \  F_{j,k} = i e_{j,k} + i e_{k,j}, \ D_{j,k} = ie_{j,j} - i e_{k,k} .
$$

By a slight abuse of language, also set 
$e_{\omega(j),\omega(k)}=e_{j,k}$. The analogous identification can be used to define $E_{\omega(j),\omega(k)},F_{\omega(j),\omega(k)},D_{\omega(j),\omega(k)}$.

Thanks to this notation we can conveniently represent the matrices corresponding 
to the 
controlled vector field (projected on $\mathcal{L}$).
A computation shows that the 
control potential in the $z$ direction, $-i \cos \theta$, projected on $\mathcal{L}$, has a matrix representation
with respect to the chosen basis 
$$
B_{3} = \sum_{m=-\l}^{\l} p_{\l,m} F_{(\l,m),(\l+1,m)}
$$
with
$$p_{\l,m}=- \sqrt{\frac{(\l+1)^{2} - m^{2}}{(2\l+1)(2\l+3)}}.$$

Similarly, we associate with the control potentials in the $x$ and $y$ directions, $-i \cos \varphi\sin \theta$ and 
$-i \sin \varphi\sin \theta$ respectively, the matrix representations
\begin{align*}
B_{1} &= \sum_{m=-\l}^{\l}( - q_{\l,m} F_{(\l,m),(\l+1,m-1)} + q_{\l,-m} F_{(\l,m),(\l+1,m+1)})\\
B_{2} &= \sum_{m=-\l}^{\l}(  q_{\l,m} E_{(\l,m),(\l+1,m-1)} +q_{\l,-m} E_{(\l,m),(\l+1,m+1)}),
\end{align*}
where
$$
q_{\l,m}= \sqrt{\frac{(\l-m+2)(\l-m+1)}{4(2\l+1)(2\l+3)}}.
$$

The matrix representation of the Schr\"odinger operator $i \Delta$  is the diagonal matrix 
$$\tilde{A}=\sum_{(j,k)\in J_\l} \tilde\alpha_{(j,k)}e_{(j,k),(j,k)}$$
where
$$
\tilde{\alpha}_{(j,k)}=-i j(j+1), \quad \mbox{ for } (j,k)\in J_\l.
$$

Now consider $A=\tilde{A} - \frac{\mathrm{tr}(\tilde{A})}{4(\l+1)} I$, in such a way that $\mathrm{tr}(A) = 0$. 
Hence, 
${A}=\sum_{(j,k)\in J_\l} \alpha_{(j,k)}e_{(j,k),(j,k)}$
where
$$
{\alpha}_{(\l,k)}=i\frac{2\l+3}{2}, \quad \mbox{ for } k=-\l,\ldots, \l,
$$
and
$$
{\alpha}_{(\l,k)}= -i\frac{2\l+1}{2},\quad \mbox{ for } k=-\l-1,\ldots, \l+1.
$$

\subsection{Useful bracket relations }\label{useful}

From the identity
\begin{equation}\label{eq:coupling}
[e_{j,k},e_{n,m}] = \delta_{kn}e_{j,m} - \delta_{jm}e_{n,k}
\end{equation}
we get the relations
$ [E_{j,k}, E_{k,n}] = E_{j,n}$, $[F_{j,k}, F_{k,n}] = -E_{j,n}$, and 
 $[E_{j,k}, F_{k,n}] = F_{j,n}$ and  
 \begin{equation}\label{eq:1231}
   [E_{j,k}, F_{j,k}] = 2D_{j,k}. 
   \end{equation}
The relations above can be interpreted following 
a ``triangle rule'': 
the bracket between an operator coupling the states $Y^{m}_{\l}$ and $Y^{n}_{k}$ and an operator coupling the states $Y^{m}_{\l}$ and $Y^{n'}_{k'}$ couples the states $Y^{n}_{k}$ and $Y^{n'}_{k'}$. 
On the other hand, the bracket is zero if two operators couple no common states.

Moreover, 
\begin{subequations}
\label{eq:parentesi}
\begin{align}
[A,E_{(\l,k),(\l+1,h)}]&=2(\l+1)F_{(\l,k),(\l+1,h)},\\
[A,F_{(\l,k),(\l+1,h)}]&=-2(\l+1)E_{(\l,k),(\l+1,h)}.
\end{align}
\end{subequations}

{
From \eqref{eq:coupling} 
we find also that
$$
[E_{(\l,m),(\l+1,m)},E_{(\l,m'),(\l+1,m'-1)}] \neq 0
$$
if and only if $m'=m$ or $m'=m+1$, with  
$$
[E_{(\l,m),(\l+1,m)},E_{(\l,m),(\l+1,m-1)}] =  E_{(\l+1,m-1),(\l+1,m)}  
$$
and 
$$
[E_{(\l,m),(\l+1,m)},E_{(\l,m+1),(\l+1,m)}]  = E_{(\l,m),(\l,m+1)}. 
$$


}

\subsection{Controllability result}

We prove the following result, which allows us to apply the abstract controllability 
criterium obtained in the previsous section. 
We obtain then Theorem~\ref{controllability} as a corollary of Theorem~\ref{metalemma}. Notice that the conclusions of Theorem~\ref{metalemma} allow us to claim more than the required approximately controllability,  since simultaneous controllability is obtained as well.

\begin{proposition}\label{prop:liemolecule}
The Lie algebra $L$ generated by $A, B_{1}, B_{2}, B_{3}$ is the whole algebra $\mathfrak{su}(4\l+4)$.
\end{proposition}

Thanks to the matrix relations obtained in Section~\ref{useful}, the proof of the proposition can be easily reduced to the proof of the following lemma.

\begin{lemma}\label{lem:1111}
The Lie algebra $L$ contains the elementary matrices
$$
E_{(\l,k),(\l+1,k+j)}\quad \mbox{ for }k=-\l,\dots,\l,\ j=-1,0,1.
$$
\end{lemma}

\noindent {\bf Proof of Lemma~\ref{lem:1111}.}
First, we want to prove that
\begin{equation}\label{eq:verticali}
\{E_{(\l,-j),(\l+1,-j)}+E_{(\l,j),(\l+1,j)} \mid j=0,\dots,\l\}\subset L. 
\end{equation}
We use the fact that 
$$
\ad^{2j+1}_{B_{3}} A = (-1)^{j} (\l+1) 2^{2j+1} \sum_{m=-\l}^{\l} p_{\l,m}^{2j+1}E_{(\l,m),(\l+1,m)}.
$$
Indeed, for $j=0$ 
\begin{align*}
[B_{3},A]  &= \sum_{\l=-m}^{m} p_{\l,m} [F_{(\l,m),(\l+1,m)},A]\\
& = 2(\l+1)   \sum_{\l=-m}^{m} p_{\l,m} E_{(\l,m),(\l+1,m)}\,,
\end{align*}
and, by induction, for $j\geq 1$,
\begin{align*}
\ad^{2j+1}_{B_{3}} A &= [B_{3},[B_{3},\ad^{2j-1}_{B_{3}} A] ]\\
& = (-1)^{j-1}(\l+1)2^{2j-1}\\
&\sum_{m=-\l}^{\l} p_{\l,m}^{2j-1} [B_{3},[B_{3},  E_{(\l,m),(\l+1,m)}]]\\
& = (-1)^{j-1}(\l+1)2^{2j-1}\sum_{m=-\l}^{\l} p_{\l,m}^{2j-1}\\
& \ \ \ [B_{3},[\sum_{h=-\l}^{\l} p_{\l,h} F_{(\l,h),(\l+1,h)},  E_{(\l,m),(\l+1,m)}]]\\
& = (-1)^{j-1}(\l+1)2^{2j-1}\\
&\sum_{m=-\l}^{\l} p_{\l,m}^{2j-1} [B_{3}, - 2 p_{\l,m} D_{(\l,m),(\l+1,m)}]\\
& = (-1)^{j} (\l+1)2^{2j} \\
&\sum_{m=-\l}^{\l} p_{\l,m}^{2j} [\sum_{h=-\l}^{\l} p_{\l,h} F_{(\l,h),(\l+1,h)}, D_{(\l,m),(\l+1,m)} ]\\
 &= (-1)^{j} (\l+1) 2^{2j+1} \sum_{m=-\l}^{\l} p_{\l,m}^{2j+1} E_{(\l,m),(\l+1,m)}.
\end{align*}
Then \eqref{eq:verticali} follows from the fact that 
$p_{\l,m} \neq p_{\l,n}$ for every $n \neq m,-m$. 

Now note that
$$
B_{2} - [A,B_{1}]/(2(\l+1)) = 2\sum_{m=-\l}^{\l} q_{\l, -m} E_{(\l,m),(\l+1,m+1)}
$$
 and
 $$
 B_{2} + [A,B_{1}]/(2(\l+1)) = 2\sum_{m=-\l}^{\l} q_{\l, m} E_{(\l,m),(\l+1,m-1)}.
 $$
 


Moreover
\begin{align*}
&{[[ \sum_{m=-\l}^{\l} q_{\l, m} E_{(\l,m),(\l+1,m-1)}, E_{(\l,0),(\l+1,0)}],E_{(\l,0),(\l+1,0)}]= }
\\
&= -q_{\l,1}[E_{(\l,0),(\l,1)},E_{(\l,0),(\l+1,0)}]\\
&\ \ \  -q_{\l,0} [E_{(\l+1,-1),(\l+1,0)},E_{(\l,0),(\l+1,0)}]\\
&=  q_{\l,1} E_{(\l,1),(\l+1,0)} + q_{\l,0}E_{(\l,0),(\l+1,-1)}.
\end{align*}
and, for $0<k\leq\l$,
\begin{align*}
&[[\sum_{j=k}^{\l} q_{\l,-j} E_{(\l,-\l),(\l+1,-\l-1)} + \\
&\ldots +q_{\l,-k+1}E_{(\l,-k+1),(\l+1,-k)}
 + q_{\l,k}E_{(\l,k),(\l+1,k-1)}+ \\
&+\ldots+ q_{\l,\l}E_{(\l,\l),(\l+1,\l-1)}, E_{(\l,-k),(\l+1,-k)}+\\
&E_{(\l,k),(\l+1,k)} ], E_{(\l,-k),(\l+1,-k)}+E_{(\l,k),(\l+1,k)} ] \\
& = -q_{\l,-k+1}[E_{(\l,-k),(\l,-k+1)},E_{(\l,-k),(\l+1,-k)}]\\
& -q_{\l,k} [E_{(\l+1,k-1),(\l+1,k)},E_{(\l,k),(\l+1,k)}]\\
 &= q_{\l,-k+1}E_{(\l,-k+1),(\l+1,-k)} + q_{\l,k} E_{(\l,k),(\l+1,k-1)}.
\end{align*}
Then we get $E_{(\l,-\l),(\l+1,-\l-1)}$, $E_{(\l,-\l+1),(\l+1,-\l)}+E_{(\l,\l),(\l+1,\l-1)}$, \ldots, $E_{(\l,0),(\l+1,-1)}+E_{(\l,1),(\l+1,0)} \in L$. Similarly we can prove that the Lie algebra $L$ contains $E_{(\l,\l),(\l+1,\l+1)}$.

Now, since  
$E_{(\l,m),(\l+1,m-1)} \in L$ and using the relation 
\begin{align*}
\mathrm{ad}_{E_{(\l,m),(\l+1,m-1)}}^2 E_{(\l,m),(\l+1,m)}+E_{(\l,-m),(\l+1,-m)}=\\
 [E_{(\l+1,m-1),(\l+1,m)},E_{(\l,m),(\l+1,m-1)}] = -E_{(\l,m),(\l+1,m)}
\end{align*}
we obtain that 
 $E_{(\l,m),(\l+1,m)}$ and $ E_{(\l,-m),(\l+1,-m)}$ belong to $L$ for every $m= -\l,\ldots,-1$

Similarly,  $E_{(\l,m),(\l+1,m)} \in L$ implies that $E_{(\l,m+1),(\l+1,m)}$ and $E_{(\l,-m),(\l+1,-m-1)}$ belong to $L$ for every $m= -\l,\ldots,-1$
\hfill$\Box$

\bibliographystyle{abbrv}
\bibliography{biblio}

\end{document}